\documentclass[12pt,a4paper]{amsart}
\usepackage{amssymb,amsmath,amsthm}
\usepackage{graphicx}

\usepackage{epsfig}

\long\def\onefigure#1#2{
\begin{figure*}[tbp]
\begin{center}
#1
\end{center}
\caption{#2}
\end{figure*}
} 

\newcommand{\lipefig}[2]  
{\onefigure{\mbox{\psfig{file=#1.eps}}}{\label{f:#1} #2} }

\begin{document}

\theoremstyle{plain}
\newtheorem{theorem}{Theorem}[section]
\newtheorem{lemma}[theorem]{Lemma}
\newtheorem{prop}[theorem]{Proposition}
\newtheorem{corollary}[theorem]{Corollary}
\newtheorem{claim}[theorem]{Claim}

\newcommand{\sgn}{\textrm{sign}}
\newcommand{\tri}{\triangle}
\newcommand{\al}{\alpha}
\newcommand{\be}{\beta}
\newcommand{\de}{\delta}
\newcommand{\De}{\Delta}
\newcommand{\ga}{\gamma}
\newcommand{\la}{\lambda}
\newcommand{\eps}{\varepsilon}
\newcommand{\N}{\mathbb{N}}
\newcommand{\R}{{\mathbb{R}^2}}
\newcommand{\z}{{\mathbb{Z}^2}}
\newcommand{\rr}{\mathbb{R}}
\newcommand{\A}{\mathrm{A}}
\newcommand{\cc}{\mathcal{C}}
\newcommand{\F}{\mathcal{F}}
\newcommand{\pp}{\mathcal{G}}
\newcommand{\dd}{\mathcal{D}}
\newcommand{\hh}{\mathcal{H}}
\newcommand{\dist}{\textrm{dist}}
\newcommand{\conv}{\textrm{conv}}
\newcommand{\aff}{\textrm{aff}}
\newcommand{\intt}{\textrm{int\;}}
\newcommand{\sign}{\textrm{sign}}

\numberwithin{equation}{section}

\title{Magic mirrors, dense diameters, Baire category}

\author{Imre B\'ar\'any, Mikl\'os Laczkovich}



\begin{abstract} An old result of Zamfirescu says that for most convex curves $C$ in the plane most points in $\R$ lie on infinitely many normals to $C$, where most is meant in Baire category sense. We strengthen this result by showing that `infinitely many' can be replaced by `contiunuum many' in the statement. We present further theorems in the same spirit.
\end{abstract}

\maketitle
\bigskip

\section{Introduction}

In a 1982 paper \cite{Za82} Tudor Zamfirescu proved  a remarkable result saying that `most mirrors are magic'. For the mathematical formulation let $\cc$ be the set of all closed convex curves in the plane $\R$. Fix some $C \in \cc$ and $z \in C$ so that the tangent line, $T(z)$, to $C$ at $z$ is unique, then so is the normal line $N(z)$ to $C$ at $z$. A point $u \in \R$ {\sl sees an image} of another point $v \in \R$ via $z$ if $u$ and $v$ and $C$ lie on the same side of $T(z)$ and the line $N(z)$ halves the angle $\angle uzv$. In particular, $u$ sees an image of itself via $z$ if $u \in N(z)$ and $u$ and $C$ are on the same side of $T(z)$.

With the Haussdorf metric $\cc$ becomes  a complete metric space. It is well-known that the normal $N(z)$ is unique at every point $z \in C$ for most convex curves $C \in \cc$ in the Baire category sense, that is, for the elements of a comeagre set of curves in $\cc$. Now the `most mirrors are magic' statement is, precisely, that for most convex curves, most points in $\R$ (again in Baire category sense) see infinitely many images of themselves. Another theorem from \cite{Za82} says that for most convex curves, most points in $\R$ see infinitely many images of any given point $v \in \R$. Zamfirescu actually proves the existence of countably many images and self-images.

The purpose of this paper is to show that most mirrors are even more magic.

\begin{theorem}\label{th:selfmirr} For most convex curves, most points in $\R$ see continuum many images of themselves.
\end{theorem}

\begin{theorem}\label{th:mirr} For most convex curves $C$ and for every point $v \in \R\setminus C$, most points in $\R$ see continuum many images of $v$.
\end{theorem}

The condition $v \notin C$ in the last theorem is used to avoid some trivial complications in the proof. The statement holds even for $v \in C$.

{\bf Remark.} Let $C^o$ denote the closed convex set whose boundary is $C$. The above definition of `$u$ sees an image of $v$ via $z\in C$' means that the mirror side of $C$ is the interior one, that is, the segment $uz$ intersects the interior of $C^o$. Theorem~\ref{th:selfmirr} does not hold when the mirror is on the other side of $C$ because every point in $\R \setminus C^o$ lies on exactly one outer normal halfline to $C$.

A statement, analogous to Theorem~\ref{th:mirr} about affine diameters was proved in \cite{BZ} in 1990, for typical $d$-dimensional convex bodies for every $d\ge 2$. The segment $[a,b]$ is an {\sl affine diameter} of $C \in \cc$ if there are distinct and parallel tangent lines to $a,b\in C$. The result in \cite{BZ} says that for most convex curves $C \in \cc$, most points on a fixed affine diameter of $C$ are contained in infinitely many affine diameters of $C$. In this case again we show the existence of continuum many diameters passing through most points in $C^o$.

\begin{theorem}\label{th:diam} For most convex curves $C \in \cc$, most points in $C^o$ lie in continuum many diameters of $C$.
\end{theorem}

Note that every point outside $C^o$ lies on the line of at most one affine diameter as any two affine diameters have a point in common. It is not hard to see, actually, that every point outside $C$ lies on a unique affine diameter.

\section{Plan of proof}

For $C \in \cc$ let $\rho(z)$ denote the radius of curvature of $C$ at $z\in C$. Let $\dd$ denote the family of all convex curves $C \in \cc$ such that
\begin{enumerate}
\item[(1)] there is a unique tangent line to $C$ at every $z \in C$,
\item[(2)] $\{z\in C: \rho(z)=0\}$ is dense in $C$,
\item[(3)] $\{z\in C: \rho(z)=\infty\}$ is dense in $C$.
\end{enumerate}

It is well-known, see for instance \cite{Za85}, that $\dd$ is comeagre in $\cc$. We are going to show that every $C \in \dd$ has the property required in Theorem~\ref{th:selfmirr}. We will need slightly different conditions for Theorems~\ref{th:mirr} and \ref{th:diam}. But the basic steps of the proofs are the same. We explain them in this section in the case of Theorem~\ref{th:selfmirr}.

Let $C \in \dd$ and define, for $z \in C$, the halfline $N^+(z)\subset N(z)$ that starts at $z$ and intersect the interior of $C^o$. Note that every $u \in \R$ lies on some $N^+(z)$: namely when the farthest point from $u$ on $C$ is $z$. Set $L(u)=\{z \in C: u \in N^+(z)\}$ and define
\[
H=\{u \in \R: L(u) \mbox{ is not perfect}\}.
\]

\begin{lemma}\label{l:borel} $H$ is a Borel set.
\end{lemma}

Write now $u=(u_1,u_2)\in \R$ and define $H^{u_2}=\{u_1\in \rr: (u_1,u_2)\in H\}$. This is just the section of $H$ on the horizontal line $\ell(u_2)=\{(x,y)\in \R: y=u_2\}$. There are two points $z \in C$ with $N(z)$ horizontal, so there are at most two exceptional values for $u_2$ where $\ell(u_2)$ coincides with some $N(z)$.

\begin{lemma}\label{l:firstcat} Apart from those exceptional values, $H^{u_2}$ is meagre.
\end{lemma}

These two lemmas imply Theorem~\ref{th:selfmirr}. Indeed, deleting the (one or two) exceptional lines from $H$ gives a Borel set $H'$. According to a theorem of Kuratowski (see \cite{Ke} page 53), if all horizontal sections of the Borel set $H'$ are meagre, then so is $H'$, and then $H$ itself is meagre. So its complement is comeagre, so $L(u)$ is perfect and non-empty for most $u \in \R$. The theorem follows now from the fact that a non-empty and perfect set has continuum many points. The proofs of Theorems~\ref{th:mirr} and \ref{th:diam} will use the same argument.

\smallskip
For the proof of Lemma~\ref{l:firstcat} we need another lemma that appeared first as Lemma 2 in \cite{Pa}. A function $g:[0,1]\to \R$ is increasing on an interval $I\subset [0,1]$ (resp. decreasing on $I$) if every $x,y \in I$ with $x \le y$ satisfy $g(x) \le g(y)$ (resp. $g(x)\ge g(y))$, and $g$ is monotone in $I$ if it is either increasing or decreasing there. For the sake of completeness we present the short proof.

\begin{lemma}\label{l:geom} Assume $g:[0,1]\to \R$ is continuous and is not monotone in any subinterval of $[0,1]$. Then the set
\[
B=\{b \in \rr: \{x:g(x)=b\} \mbox{ is not perfect}\}
\]
is meagre.
\end{lemma}

\medskip
{\bf Proof} of Lemma~\ref{l:geom}. For each $b \in B$ the level set $\{x:g(x)=b\}$ has an isolated point, and so there is an open  interval $I_b \subset [0,1]$ with rational endpoints in which $g(x)=b$ has a unique solution. For a given rational interval $(p,q)$ define 
\[
B(p,q)=\{b\in B: I_b=(p,q)\}.
\]
If every $B(p,q)$ is nowhere dense, then we are done since $B$, as a countable union of nowhere dense sets, is meagre. If some $B(p,q)$ is not nowhere dense, then there is a non-empty open interval $I$ in which $B(p,q)$ is dense. The line $y=b$, for a dense subset of $I$, intersects the graph of $g$ restricted to $(p,q)$ in a single point. This implies easily that $g$ is strictly monotone in a subinterval $(p,q)$, contrary to our assumption.\hfill$\Box$

\section{Proof of the lemmas}

Fix $C \in \dd$ and let $z(\al)$ denote the point $z \in C$ where the halfline $N^+(z)$ spans angle $\al\in [0,2\pi)$ with a fixed unit vector in $\R$. This is a parametrization of $C$ with $\al \in [0,2\pi]$ and $z(0)=z(2\pi)$. We write $C_{\al,\be}$ for the arc $\{z(\ga): \al < \ga < \be\}$ when $0\le \al < \be \le 2\pi$, and the definition is extended, naturally, to the case when $\al < 2\pi < \be$. We always assume that $\al, \be$ are rational and $\be -\al$ is small, smaller than $0.1$, say.

\medskip
{\bf Proof} of Lemma~\ref{l:borel} . Note first that the set
\[
K=\{(u,z)\in \R \times C: u \in N^+(z)\}
\]
is closed. Further, $L(u)$ is not perfect if and only if there is a short arc $C_{\al,\be}$ such that $u \in N^+(z)$ for a unique $z \in C_{\al,\be}$. Thus
\[
H = \bigcup_{\mbox{ all }C_{\al,\be}} \{u \in \R: u \in N^+(z) \mbox{ for a unique } z \in C_{\al,\be}\}.
\]
Let $p: K \to \R$ be  the projection $p(u,z)=u$. Let $P_{\al,\be}$ be the set of points $u \in \R$ such that there are more than one $z \in C_{\al,\be}$ with $u\in N^+(z)$. Then
\[
P_{\al,\be}=\bigcup_{\ga} p(K\cap (\R \times  C_{\al,\ga}))\cap p(K\cap (\R \times  C_{\ga,\be}))
\]
where the union is taken over all rational $\ga$ with $\al<\ga < \be$. Consequently
\[
H=\bigcup_{\mbox{ all }C_{\al,\be}} p(K\cap (\R \times  C_{\al,\be})) \setminus P_{\al,\be}.
\]
Since $p(K\cap (\R \times  C_{\al,\be}))$ is $F_{\sigma}$ for every $\al < \be$, it follows that $H$ is indeed Borel. \hfill$\Box$

\medskip
{\bf Proof} of Lemma~\ref{l:firstcat}. The set $z\in C$ where $N^+(z)$ intersects $\ell(u_2)$ in a single point consists of one or two open subarcs of $C$, as one can check easily. Let $C_1$ be such an arc. It suffices to see that
\[
E=H^{u_2}\cap \{u_1\in \rr: (u_1,u_2)=\ell(u_2)\cap N(z) \mbox{ for some }z \in C_1\}
\]
is meagre, as $H^{u_2}$ either coincides with this set, or is the union of two such sets.

We may assume that $C_1$ is the graph of a convex function $F: J \to \rr$ and $u_2>F(x)$ on $J$ where $J$ is an open interval.
(This position can be reached after a suitable reflection about a horizontal line.) With this notation, $E$ is the set of real numbers $u_1 \in \rr$ such that the set of points $x\in J$ for which $(u_1,u_2)\in N^+(x,F(x))$ is not perfect.

Then $F'(x)=f(x)$ is continuous and increasing on $J$. Each $z \in C_1$ is a point $(x,F(x))$ on the graph of $F$. As $\rho(z) =(1+f(x))^{3/2}/f'(x)$, $f'$ equals zero resp. infinity on a dense set in $J$. The normal $N(z)$ to $z=(x,F(x))$ has equation $(u_2-F(x))f(x)=x-u_1$,
as one checks readily. With $g(x)=(u_2-F(x))f(x)-x$, $g'(x)=-f(x)^2+(u_2-F(x))f'(x)-1$ and so on a dense set in $J$ the value of $g'(x)$ is positive, and on another dense set in $J$ it is negative. So $g$ is not monotone in any subinterval of $J$. Lemma~\ref{l:geom} implies now that $E$ is meagre. \hfill$\Box$

\section{Proof of Theorem~\ref{th:mirr}}

It is known \cite{Za85} that for most $C \in \dd$ there is a dense set $E\subset C$ such that at each point $z \in E$ the lower curvatures of radii in both directions $\rho_i^+(z),\rho_i^-(z)$ vanish and the upper curvatures of radii in both directions $\rho_s^+(z),\rho_s^-(z)$ are infinite. We let $\dd_1$ denote the set of all $C \in \dd$ possessing such a dense set $E$. We are going to show that for each $C \in \dd_1$, most points see continuum many images of any given point $v \in \R$, $v\notin C$.

For $z \in C$ we define the line $R(z)$ as the reflected copy (with respect to $N(z)$) of the line through $v$ and $z$. Note that $R(z)$ depends continuously from $z$. Here we need $v \notin C$.

If $u$ sees an image of $v$ via $z$, then $u \in R(z)$. More precisely, $u$ sees an image of $v$ via $z$ iff $u,v$ and $C$ are on the same side of $T(z)$ and $u \in R(z)$. Let $R^+(z)\subset R(z)$ be the halfline that starts at $z$ and is on the same side of $T(z)$ as $C$. Also, $R^+(z)$ is well defined for all $z \in C$.

As before, $\ell(u_2)$ is the horizontal line in $\R$ whose points have second coordinate equal to $u_2$. Define, for fixed $u_2 \in \rr$, $H^{u_2}=\{u_1 \in \rr: (u_1,u_2)\in H\}$. This is the same as the set of first coordinates of all $u \in H\cap \ell(u_2)$.

In the generic case $R(z)$ is not horizontal and so $R(z)\cap \ell(u_2)$ is a single point. But we have to deal with non-generic situations, that is, when $R(z)$ is horizontal and so coincides with $\ell(u_2)$ for some $u_2\in \rr$. Define $Z=\{z\in C: R(z) \mbox{ is horizontal}\}$ and $U_2=\{u_2 \in \rr: \ell(u_2)=R(z) \mbox{ for some } z \in Z\}$. Both $Z$ and $U_2$ are closed sets and there is a one-to-one correspondence between them given by $z \leftrightarrow u_2$ iff $R(z)=\ell(u_2)$.

From now on we assume that $Z$ is nowhere dense. We will justify this assumption at the end of the proof. Then $U_2$ is also nowhere dense. $C\setminus Z$ is open in $C$ and so its connected components $C_1,C_2,\ldots$ are open arcs in $C$, and there are at most countably many of them.

This time we define $L(u,C_i)$ as the set of $z \in C_i$ via which $u$ sees an image of $v$. Formally, $L(u,C_i)=\{z \in C_i: u \in R^+(z)\}$, and define again, for fixed $u_2 \in \rr$,
\[
H_i^{u_2}=\{u_1\in \rr: L((u_1,u_2),C_i) \mbox{ is not perfect}\}.
\]

A very similar proof shows that $H_i^{u_2}$ is Borel. We omit the details, but mention that the condition $v \notin C$ is needed to show that the corresponding $K=\{(u,z):\dots\}$ is closed. 

\medskip
\begin{lemma}\label{l:firstcat2} For $u_2 \notin U_2$ the set $H_i^{u_2}$ is meagre.
\end{lemma}

\medskip
{\bf Proof.} With every $u_1 \in H_i^{u_2}$ we associate a (rational) open arc $C_{\al,\be}$ of $C_i$ such that $u=(u_1,u_2) \in R(z)$ for a unique $z \in C_{\al,\be}$, namely for $z_u$. If the set of $u\in H_i^{u_2}$ that are associated with $C_{\al,\be}$ is nowhere dense for every rational arc $C_{\al,\be}$, then we are done as $H_i^{u_2}$ is the countable union of nowhere dense sets. So suppose that it is not nowhere dense for some $C_{\al,\be}$. Then there is an open interval $I \in \rr$ such that $H_i^{u_2}$ is dense in $I$.

Choose two distinct points $w^-,w^+$ from $I \cap H_i^{u_2}$. Then $z_{(w^-,u_2)}$ and $z_{(w^+,u_2)}$ are distinct points and so they are the endpoints of an open subarc $C_{\ga,\de}$ of $C_{\al,\be}$. Define the map $h:C_{\ga,\de} \to I$ by $h(z)=u_1$ when $(u_1,u_2)=\ell(u_2)\cap R(z)$; $h$ is clearly continuous. It is also monotone because its inverse is well-defined on a dense subset $I$.

We show next that this is impossible. Choose $z_0 \in C_{\ga,\de} \cap E$ (recall that $E$ is dense in $C$).

\begin{figure}
\centering
\includegraphics{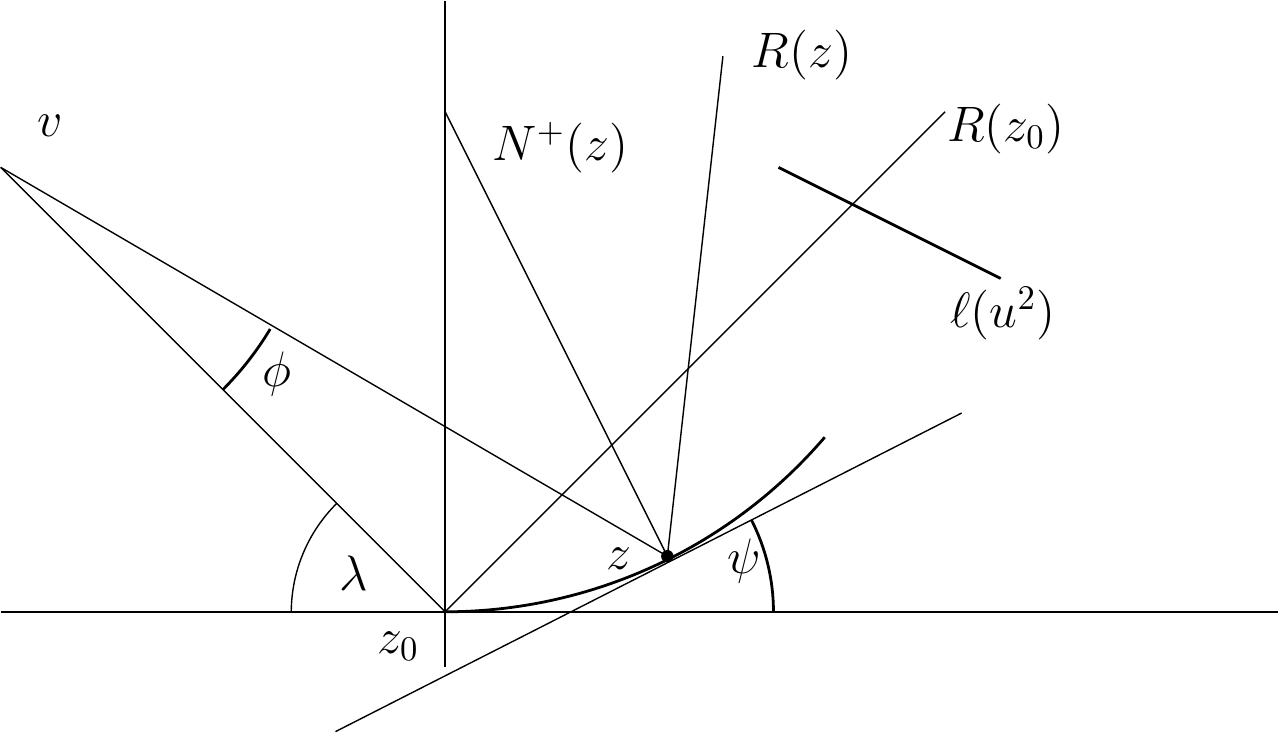}
\caption{Theorem~\ref{th:mirr}}
\label{fig:delta}
\end{figure}

We fix a new coordinate system in $\R$: the origin coincides with $z_0$, the $x$ axis with $T(z_0)$, the tangent line to $C$ at $z_0$, and the $y$ axis is $N(z_0)$; see the figure. We assume w.l.o.g. that $v_1<0$ and $v_2>0$ where $v=(v_1,v_2)$. A subarc of $C_{\ga,\de}$ is the graph of a non-negative convex function $F:[0,\De)\to \rr$ such that $F(0)=0$ and $z=z(x)=(x,F(x))$ and $f(x)=F'(x)$ is an increasing function with $f(0)=0$. If the lines $R(z(x))$ and $R(z(0))$ intersect, then they intersect in a single point whose $y$ component is denoted by $y(x)$.

\begin{claim}\label{cl:eps} For every $\eps>0$ there are $x_1,x_2 \in (0,\eps)$ so that $y(x_1)<0$ and $0<y(x_2)<\eps$.
\end{claim}

{\bf Proof.} We use the notation of the figure. The sine theorem for the triangle with vertices $v,0,z(x)$ implies that $\phi(x)=x \sin \la/|v|(1+o(1))$ where $o(1)$ is understood when $x \to 0$. The slope of the line $R(z(x))$ is $\tan (\la-\phi+2\psi)$, and
\[
\tan \psi(x)=f(x)=x\cdot \frac {f(x)-0}{x-0}.
\]
The liminf and limsup of the last fraction (when $x\to 0$) is the curvature $\rho^+_i(z_0)=0$ and $\rho^+_s(z_0)=\infty$ of $C$ at $z_0$ as $z_0 \in E$. Consequently for every integer $n>1$ there is $x\in (0,1/n)$ with $\tan \psi(x)< x/n$ and also with $\tan \psi(x)>nx$. Then there is $x_1<1/n$ such that $\la/2 < \la -\phi(x_1) +2\psi(x_1)< \la$ which implies, after a simple checking, that $y(x_1)<0$. Also, there is $x_2<1/n$ such that $\la -\phi(x_2) +2\psi(x_2)> \la +nx_2/2$. A direct computation shows then that $0<y(x_2)<\eps$ if $n$ is chosen large enough.\hfill$\Box$

\medskip
We return to the proof of Lemma~\ref{l:firstcat2}. The claim shows that there are $x_1,x_2,x_3 \in (0,\De)$ with $x_1<x_2<x_3$ such that the line $R(z(x_1))$ and $R(z(x_3))$ strictly separate the origin and the point $R(z_0)\cap \ell(u_2)$ while $R(z(x_2))$ does not. Writing $z_i=z(x_i), i=1,2,3$ this implies that $z_2$ is between $z_1$ and $z_3$ while $h(z_2)$ is not on the segment $(h(z_1),h(z_3))$. So $h$ is not monotone.\hfill$\Box$

\medskip
It is evident that $U_2$, and consequently $U$, is closed and nowhere dense, so $U$ is meagre. The lemma implies, by Kuratowski's theorem, that $H_i\setminus U$ is meagre. It follows that $H_i$ is meagre and then so is $H=\bigcup_i H_i$. Thus every point in the complement of $H$ sees an image of $v$ via a perfect set in $C$, except possibly for the points of the meagre set $U$. This perfect set is nonempty, because every point sees an image of $v$ via some $z \in C$ (for instance by Zamfirescu's result \cite[Theorem 1]{Za82}). So most points see continuum many images of $v$.

Finally we justify the assumption that $Z$ is nowhere dense. This is done by choosing the horizontal direction (which is at our liberty) suitably. So for a given direction $(\cos \theta,\sin \theta)$ write $Z(\theta)$ for the set of $z\in C$ such that $R(z)$ is parallel with this direction. Every $Z(\theta)$ is closed and so there is one (actually, many) among them that contains no non-empty open arc of $C$. Choose the corresponding $\theta$ for the horizontal direction, then $Z=Z(\theta)$ is nowhere dense.
\hfill$\Box$

\section{Proof of Theorem~\ref{th:diam}}

Write $\cc_1$ for the set of all convex curves $C$ that have a unique tangent at every $z \in C$.
Assume $C \in \cc_1$ and use the parametrization $z:[0,2\pi)\to C$ as before. For $z\in C$ with $z=z(\al)$ let $z^*\in C$ be the opposite point, that is $z^*=z(\al+\pi)$. It is evident that $z^{**}=z$. Further, $[z,z^*]$ is always an affine diameter of $C$ and all affine diameters of $C$ are of this form. We need a geometric lemma.

\begin{lemma}\label{l:arc} Most convex curves $C \in \cc_1$ have the following property: for every $\eps> 0$  every subarc $C_0$ of $C$ contains points $x,y$ such that
\[
\frac {|x-y|}{|x^*-y^*|} < \eps.
\]
\end{lemma}

The lemma follows from a result in \cite{AZ}, we give a separate proof in the next section. From  now on we assume that $C \in \cc_1$ has the property in the lemma.

We use again the same proof scheme: for $u \in C^o$ define $L(u)=\{z\in C: u\in [z,z^*]\}$; this set is nonempty as one can check easily that every point $u \in C^o$ lies on at least one affine diameter. (This holds for every convex curve, not only for the ones in $\cc_1$.) We set next $H=\{u \in C^o: L(u) \mbox{ is not perfect}\}$, and, for fixed $u_2 \in \R$, $H^{u_2}=H \cap \ell(u_2)$. The same proof as in Section 3 shows that $H$ is Borel. We claim that $H$ is meagre which implies Theorem~\ref{th:diam}.

$C$ has a horizontal affine diameter and we assume w.l.o.g. that it lies on the line $\ell(0)$.
To see that $H$ is meagre it suffices to show (by Kuratowski's theorem) that $H^{u_2}$ is meagre as a subset of $\ell(u_2)$ for $u_2\ne 0$. We only consider $u_2 \in \rr$, $u_2\ne 0$ with $\ell(u_2) \cap C \ne \emptyset$. With each $u \in H^{u_2}$ we associate an isolated point $z_u \in C$ and a short rational arc $C_{\al,\be}$ such that
$z_u$ is the unique $z \in C_{\al,\be}$ with $u \in [z,z^*]$. We are done if, for each short rational arc $C_{\al,\be}$, the set of $u \in H^{u_2}$ that are associated with $C_{\al,\be}$ is nowhere dense. So suppose that this fails for some $C_{\al,\be}$. Then there is an open interval $I \subset \ell(u_2)$ on which $H^{u_2}$ is dense. Choose distinct points $u^-$ and $u^+$ from $I \cap H^{u_2}$ and let $z^-,z^+$ be the corresponding isolated points on $C_{\al,\be}$. We suppose (by symmetry) that $C_{\al,\be}$ is below the line $\ell(u_2)$.

From now on we consider the subarc $C_0 \subset C_{\al,\be}$ whose endpoints are $z^-$ and $z^+$ and its opposite arc $C_0^*$.
We note here that the map $z\to z^*$ is order preserving on $C_0$, meaning that if $v\in C_0$ is between $v_1,v_2 \in C_0$, then $v^*$ lies between $v_1^*$ and $v_2^*$ on $C_0^*$.

Define a map $m:C_0 \to \ell(u_2)$ via $m(z)=\ell(u_2)\cap [z,z^*]$; $m$ is continuous. It is one-to-one on a dense subset of $C_0$ which implies that $m$ is order-preserving in the sense that if $v\in C_0$ is between $v_1,v_2 \in C_0$, then $m(v)$ lies between $m(v_1)$ and $m(v_2)$ on $\ell(u_2)$. We show that this is impossible.

Using Lemma~\ref{l:arc} choose two points $v_1,v_2$ on $C_0$ very close to each other so that $|v_1-v_2|$ is much shorter than $|v_1^*-v_2^*|$. Then the segment $[v_1,v_2]$ is almost parallel with $[v_1^*,v_2^*]$, and the diameters $[v_1,v_1^*]$ and $[v_2,v_2^*]$ intersect in a point very close to $[v_1,v_2]$, so this point is below $\ell(u_2)$. Now apply Lemma~\ref{l:arc} on the arc between $v_1^*$ and $v_2^*$. We get points $w_1$ and $w_2$ very close to each other on this arc so that $|w_1-w_2|$ is much shorter than $|w_1^*-w_2^*|$. This time the diameters $[w_1,w_1^*]$ and $[w_2,w_2^*]$ intersect above $\ell(u_2)$. We assume (by choosing the names $w_1,w_2$ properly) that $v_1^*,w_1,w_2,v_2^*$ come in this order on $C_0^*$ and so $v_1,w_1^*,w_2^*,v_2$ come in this order on $C_0$. The order of their $m$-images on $\ell(u_2)$ is $m(v_1),m(w_2^*),m(w_1^*),m(v_2)$. Thus indeed, $m$ is not order preserving.\hfill$\Box$

\section{Proof of Lemma~\ref{l:arc}}

Given $C \in \cc_1$ define $A_{k,n}$ as the short arc between $z_k=z(2\pi k/2n)$ and $z_{k+1}=z(2\pi (k+1)/2n)$ where $k=0,1,\ldots,2n-1$. For positive integers $n,m$ let $\F_{n,m}$ be the set of all $C \in \cc_1$ for which there is $A_{k,n}$ such that for all $x,y \in A_{k,n}$ ($x\ne y$)
\[
\frac{|x-y|}{|x^*-y^*|}\ge \frac 1m.
\]

It is easy to see that $\F_{n,m}$ is closed in $\cc_1$, we omit the details. We show next that it is nowhere dense.

Fix a $C \in \cc_1$ and $\eps>0$ and let $U(C)$ denote the $\eps$-neighbourhood of $C$. We construct another convex curve $\Gamma \in \cc_1$ that is contained in $U(C)$ but is not an element of $\F_{n,m}$. Fix $k\in \{0,1,\ldots,n-1\}$ and consider a fixed arc $A_{k,n}$ and its opposite arc $A^*_{k,n}=A_{k+n,n}$. Let $T_k$ be the tangent line to $C$ at $z((k+\frac 12)\pi/n)$ and $T_k^*$ be the parallel tangent line at $z((k+n+\frac 12)\pi/n)$. Translate $T_k^*$ a little so that the translated copy intersects  $C$ in two points $x_1,y_1$ and the segment $[x_1,y_1]$ lies in $U(C)$ and is much shorter than $[z_{k+n},z_{k+1+n}]$. Similarly translate $T_k$ a little so that the translated copy intersects $C$ in $x_2,y_2$ and $[x_2,y_2]$ lies in $U(C)$, and is much shorter than $[z_k,z_{k+1}]$ and, most importantly, it is much shorter than $[x_1,y_1]$, namely, $m|x_2-y_2| <|x_1-y_1|$. This is clearly possible.

Now we choose points $w_1$ resp. $w_2$ from the caps cut off from $C^o$ by the segment $[x_1,y_1]$ and $[x_2,y_2]$ so that, for $i=1,2$, the triangles $\tri_i=\conv \{x_i,y_i,w_i\}$ are homothetic. This is possible again. Note that $[x_1,w_1]$ and $[x_2,w_2]$ are parallel, and so are $[y_1,w_1]$ and $[y_2,w_2]$.

The next target is construct a convex curve $\Gamma_k$ from $z_k$ to $z_{k+1}$ going through $x_2$ and $y_2$ that lies in $U(C)$, has a unique tangent at every point, and this tangent coincides with the line through $x_2,w_2$ at $x_2$ and with the line through $y_2,w_2$ at $y_2$. Also, an analogous curve $\Gamma_{k+n}$ is needed from $z_{k+n}$ to $z_{k+1+n}$.

This is quite easy. The unique parabola arc connecting $x_2$ to $y_2$ within $\tri_2$ that touches the sides $[x_2,w_2]$ at $x_2$ and $[w_2,y_2]$ at $y_2$ is the middle piece of $\Gamma_k$. To connect this arc by a convex curve to $z_k$ (say) within $U(C)$ choose a point $w\in C$ on the arc between $z_k$ and $y_2$ so close to $y_2$ that the triangle $\tri$ delimited by $T(z)$, the line through $y_2,w_2$, and the segment $[y_2,z]$ lies in $U(C)$. The analogous parabola arc in $\tri$ gives the next piece of $\Gamma_k$, and then add to this piece the subarc of $C$ between $w$ and $z_k$. The middle piece of $\Gamma_k$ is continued to $z_{k+1}$ the same way.

The convex curve $\Gamma_{k+n}$ connecting $z_{k+n}$ to $z_{k+1+n}$ is constructed the same way. Note that the tangents to $\Gamma_k$ at $x_2$ (resp. $y_2$) are parallel with the tangents to $\Gamma_{k+n}$ at $x_1$ (and $y_1$).

The curves $\Gamma_k$ for $k=0,\ldots,2n-1$ together form a convex curve $\Gamma \in C_1$. It has parallel tangents at $x_1\in \Gamma_{k+n}$ and  $x_2\in \Gamma_k$, and also at $y_1$ and $y_2$. Thus $[x_1,x_2]$ and $[y_1,y_2]$ are affine diameters of $\Gamma$ and $m|x_1-y_1|<|x_2-y_2|$. As this holds for every $k$,  $\Gamma \notin \F_{n,m}$. Thus $\F_{n,m}$ is indeed nowhere dense.

It follows that $\cc_2=\cc_1 \setminus \bigcup_{n,m}\F_{n,m}$ is comeagre in $\cc_1$. We show next that every $C \in \cc_2$ satisfies the requirement of the lemma. So we are given $\eps>0$ and a short subarc $C_0$ of $C$. Take a positive integer $m$ with $1/m<\eps$. For a suitably large $n$, $C_0$ contains an arc of the form $A_{k,n}$. As $C \notin \F_{n,m}$, there are distinct points $x,y \in A_{k,n}$ with
\[
\frac{|x-y|}{|x^*-y^*|}\le \frac 1m< \eps.
\]
This finishes the proof. \hfill$\Box$

\bigskip
{\bf Acknowledgements.} Research of the first author was partially supported by ERC Advanced Research Grant 267165 (DISCONV), and by Hungarian National Research Grant K 83767. The second author was partially supported by the Hungarian National Foundation for Scientific Research Grant K 104178.

\end{document}